\documentclass{article}
\title{The structure of ${\cal{AK}}_2$-manifolds.
}
\author{Paul-Andi Nagy}
\date{\today}
\usepackage{amsfonts,amssymb}
\usepackage{minitoc_href}

\newtheorem{teo}{Theorem}[section]
\newtheorem{lema}{Lemma}[section]
\newtheorem{pro}{Proposition}[section]
\newtheorem{defi}{Definition}[section]
\newtheorem{rema}{Remark}[section]

\newtheorem{nr}{}[section]
\begin{document}
\maketitle
\abstract{ \normalsize We study special almost K\"ahler manifolds whose curvature tensor satisfies the second 
curvature condition of Gray. It is shown that for such manifolds, the torsion of the first canonical Hermitian is 
parallel. This enables us to show that every ${\cal{AK}}_2$-manifold has parallel torsion. Some applications  of this 
result, concerning the existence of orthogonal almost K\"ahler structures on spaces of constant curvature, are given. 
\large
\tableofcontents
\section{Introduction}
An  almost K\"ahler manifold (shortly ${\cal{AK}}$) is a Riemannian manifold $(M^{2n},g)$, together with a compatible almost complex structure 
$J$, such that the K\"ahler form $\omega=g(J \cdot, \cdot)$ is closed. Hence, almost K\"ahler geometry is nothing else that symplectic 
geometry with a prefered metric and complex structure. Since symplectic manifolds often arise in this way, it is rather 
natural to ask under which conditions on the metric we get integrability of the almost complex $J$. In this direction, 
a famous conjecture of S. I. Golberg asserts that every compact, Einstein, almost K\"ahler manifold is, in fact, K\"ahler. At our present 
knowledge, this conjecture is still open. Nevertheless, they are a certain number of partial results, supporting this 
conjecture. First of all, K. Sekigawa proved \cite{Seki1} that 
the Goldberg conjecture is true when the scalar curvature is positive.
We have to note that the Golberg conjecture is definitively not true with the compacity assumption removed. In fact, there are 
Hermitian symmetric spaces of non-compact type of any complex dimension $n \ge 3$ admitting almost K\"ahler structures commuting with the invariant 
K\"ahler one \cite{Apo2}. Not that, at the opposite, the real hyperbolic space of dimension at least $4$ do not admits, even locally, orthogonal almost K\"ahler structures 
\cite{Olszak, arm2}. In dimension $4$, examples of local Ricci flat almost K\"ahler metrics are constructed in 
\cite{Apo, arm2,Nur} . 
In the same paper, a potential source of compact almost K\"ahler, Einstein manifolds is considered, 
namely those compact K\"ahler manifolds whose Ricci tensor admits two distinct, constant eigenvalues; integrability 
is proven under certain positivity conditions. 
The rest of known results, most of them 
enforcing or replacing the Einstein condition with some other natural curvature assumption are mainly 
in dimension $4$. To cite only a few of them, we mention the beautiful series of papers \cite{Apo1, Apo3, Apo4} giving 
a complete local and global classfication of almost K\"ahler manifolds of $4$ dimensions satisfying the second and third 
Gray condition on the Riemannian curvature tensor. Other recent results, again in $4$-dimensions, are concerned with the study of local obstructions 
to the existence of Einstein metrics \cite{arm2}, $\star$-Einstein metrics \cite{Seki2}, etc.
\par
The aim of this paper is to close the circle of ideas from our previous paper \cite{Nagy} and to show coincidence between the class 
of ${\cal{AK}}_2$-manifolds- that is almost K\"ahler manifolds whose curvature tensor satisfies the second Gray's condition- and almost K\"ahler manifolds with 
parallel torsion. Here, both torsion and parallelism are with respect ot the first canonical Hermitian connection. 
More precisely, let us recall that it was proven in \cite{Nagy} that in an open dense set every ${\cal{AK}}_2$ locally splits a the 
Riemannian product of an almost K\"ahler manifolds with parallel torsion and a special ${\cal{AK}}_2$-manifold (see section 
3 for the definition). Our main present result 
is to give a classification of special ${\cal{AK}}_2$-manifolds. 
\begin{teo}
Let $(M^{2n},g,J)$ be a special ${\cal{AK}}_2$-manifold. Then $M$ has parallel torsion and $(M,g,J)$ is in fact a locally $3$-symmetric 
space.
\end{teo}
Therefore we get the precise description, in terms of the torsion of almost K\"ahler manifolds whose curvature tensor satisfies the 
second Gray condition. 
\begin{teo}
Let $(M^{2n},g,J)$ be an almost K\"ahler manifold in the class ${\cal{AK}}_2$. Then 
the torsion of $(g,J)$ is parallel with respect to the first canonical Hermitian connection. 
\end{teo}
As a corollary we obtain that there are no Einstein manifolds in the class ${\cal{AK}}_2$. This follows from \cite{Nagy} as a simple application of 
Sekigawa's formula combined with the paralleism of the torsion. 
Let us note that our actula classification of ${\cal{AK}}_2$-manifolds is of course not a complete one. What it remains 
to be done is to classify almost K\"ahler manifolds with parallel torsion with respect to the first 
canonical Hermitian connection. \par
Our present classification enables some applications relating to the non-existence of almost K\"ahler structures in 
hyperbolic geometry.
\begin{teo}
(i) Let $(M,g)$ be a manifold of constant sectional curvature. If there exists an almost complex structure $J$ such that 
$(g,J)$ is almost K\"ahler, then $g$ is a flat metric. \\
(ii) Let $(M^{2n},g,J)$ be a K\"ahler manifold of constant negative holomorphic sectional curvature. If $I$ is an 
almost complex structure commuting with $I$ such that $(g,I)$ is almost K\"ahler, then $I$ has to be in fact a 
K\"ahler structure.
\end{teo}
The result in (i) is not new, with exception of its proof. In dimension beyond $8$ it was proven in \cite{Olszak}, and 
in dimension $4$ and $6$ in \cite{arm2}.
Note that in dimension $4$ the result of (ii) was already obtained by J. Armstrong, as a consequence of his classification 
of almost K\"ahler, Einstein, $4$-manifold of class ${\cal{AK}}_3$. 
In dimension $6$ and beyond, the result is, at our knowledge, new. Also, note that, by results in \cite{Apo1}, that they 
are Hermitian symmetric spaces of 
non-compact type, admitting a reversing strictly almost K\"ahler structure. Other symmetric spaces can also 
support strictly almost K\"ahler structures \cite{Oguro}. \par
The paper is organised as follows. In section 2 we review some elementary facts from almost K\"ahler geometry and 
also a few results from our previous paper \cite{Nagy}, to be used in the subsequent. Section 3 is devoted to the 
investigation of the curvature tensor of special ${\cal{AK}}_2$ manifolds. Some important technical tools are developed and at 
the end of the section it is proved that any special ${\cal{AK}}_2$-manifold with parallel torsion has to be 
locally $3$-symmetric. The fourth section contains the proof of the theorem 1.1, whose basic ingredient is the observation 
that a partial holonomy reduction (with respect to the canonical Hermitian connection) extends in a canonical way 
to a global one. 
\section{Preliminaries}
Let us consider an almost Hermitian manifold $(M^{2n},g,J)$, that is a Riemannian manifold 
endowed with a compatible complex structure. We denote by $\nabla$ the Levi-Civita 
connection of the Riemannian metric $g$. Consider now the tensor $\nabla J$, the first derivative 
of the almost complex structure and recall that for all $X$ in $TM$ we have that 
$\nabla_XJ$ is a skew-symmetric (with 
respect to $g$) endomorphism of $TM$, which anticommutes with $J$. The tensor 
$\nabla J$  can be used to distinguish various classes of almost 
Hermitian manifolds. For example, $(M^{2n},g,J)$ is quasi-K\"ahler iff 
$$\nabla_{JX}J=-J\nabla_XJ $$
for all $X$ in $TM$. If $\omega=g(J \cdot, \cdot)$ denotes the K\"ahler form of the 
almost Hermitian structure $(g,J)$, we have an almost K\"ahler structure (${\cal{AK}}$ for 
short), iff $d \omega=0$. We also recall the well known fact that almost K\"ahler manifolds 
are always quasi-K\"ahler. \par
The almost complex structure $J$ defines a Hermitian structure if it is integrable, that is 
the Nijenhuis tensor $N_J$ defined by 
$$ N_J(X,Y)=[JX,JY]-[X,Y]-J[X,JY]-J[JX,Y]$$
for all vector fields $X$ and $Y$ on $M$ identically vanishes. This is also equivalent to 
$$ \nabla_{JX}J=J\nabla_XJ$$
whenever $X$ is in $TM$. Therefore, an almost K\"ahler manifold which is also 
Hermitian must be K\"ahler. \par
In the rest of this section $(M^{2n},g,J)$ will be an almost K\"ahler manifold. 
We begin to recall some basic facts about the various notions 
of Ricci tensors.  \par
Let $Ric$ be the Ricci tensor of the Riemannian metric $g$. We denote by $Ric^{\prime}$ and 
$Ric^{\prime \prime}$ the $J$-invariant resp. the $J$-anti-invariant part of the tensor $Ric$. Then the 
Ricci form is defined by 
$$ \rho=<Ric^{\prime}J \cdot, \cdot>.$$
We define the $\star$-Ricci form by 
$$ \rho^{\star}=\frac{1}{2} \sum \limits_{i=1}^{2n}R(e_i, Je_i)$$
where $\{ e_i, 1 \le i \le 2n \}$ is any local orthonormal basis in $TM$. Note that $\rho^{\star}$ is not, in general,
$J$-invariant. The $\star$-Ricci form is related 
to the Ricci form by 
\begin{nr} \hfill 
$ \rho^{\star}-\rho=\frac{1}{2} \nabla^{\star} \nabla \omega. \hfill $
\end{nr}
The (classical) proof of this fact consists in using the Weitzenb\"ock formula for the harmonic 
$2$-form $\omega$. Taking the scalar product with $\omega$ we obtain : 
$$ s^{\star}-s=\frac{1}{2}\vert \nabla J\vert^2$$
where the $\star$-scalar curvature is defined by $s^{\star}=2<R(\omega), \omega>$. \par
A basic object in almost K\"ahler geometry is the first canonical Hermitian connection, defined by : 
$$ \overline{\nabla}_XY=\nabla_XY+\eta_XY$$
for all vector fields $X$ and $Y$ on $M$. Here, the tensor $\eta$ is given by $\eta_XY=\frac{1}{2}(\nabla_XJ)JY$. Then $\overline{\nabla}$ is a metric 
Hermitian connection, that is it respects the metric and the almost complex structure. The torsion of $\overline{\nabla}$ is defined by 
$T_XY=\eta_XY-\eta_YX$. Now, it is worthly to note that the almost K\"ahler condition (i.e. $d\omega=0$) is equivalent to 
$$ <T_XY ,Z>=-<\eta_ZX,Y>$$
for all $X,Y$ and $Z$ in $TM$. In the subsequent we will refer simply to the tensor $T$ as the torsion of the almost K\"ahler manifold $(M^{2n},g,J)$.
Our main object of study in this paper is the class of almost K\"ahler manifolds introduced by the following definition. 
\begin{defi}
An almost K\"ahler manifold $(M^{2n},g,J)$ belongs to the class ${\cal{AK}}_2$ iff 
its Riemannian curvature tensor satisfies the identity : 
$$R(X,Y,Z,U)-R(JX,JY,Z,U)=R(JX,Y,JZ,U)+R(JX,Y,Z,JU)  $$
for all $X,Y,Z,U$ in $TM$.
\end{defi}
Manifolds within this class have a simple caracterization in terms of the torsion of the canonical Hermitian connection. 
\begin{pro} \cite{Nagy} 
An almost K\"ahler manifold $(M^{2n},g,J)$ belongs to the class ${\cal{AK}}_2$ iff 
$$ (\overline{\nabla}_X\eta)(Y,Z)= (\overline{\nabla}_Y\eta)(X,Z)$$ 
for all vector fields $X,Y$ and $Z$ on $M$.
\end{pro}
A number of usefull properties can be derived from the previous caracterization. At first let us define the K\"ahler nullity of the almost K\"ahler structure $(g,J)$ 
to be $H=\{ v : \eta_v=0\}$. The orthogonal complement of $H$ in $TM$ will be denoted by 
${\cal{V}}$. Since our study is purely local we can assume without loss of generality that $H$ (and hence 
${\cal{V}}$) have constant rank (this happens anyway in any connected component of a dense open subset of $M$). Then every closed property proved locally 
will extend to the whole manifold. 
\begin{pro}\cite{Nagy} We have : \\
(i) both distributions ${\cal{V}}$ and $H$ are integrable. \\
(ii) $\overline{\nabla}_V\eta=0$ for all $V$ in ${\cal{V}}$. \\
(iii) for any $V,W$ in ${\cal{V}}$ and $X$ in $H$ we have that $\overline{\nabla}_VW$ and $\overline{\nabla}_VX$ belong to ${\cal{V}}$ and $H$ respectively. 
\end{pro}
Let us denote now by $\overline{R}$ the curvature tensor of connection $\overline{\nabla}$. It has the following symmetry property : 
\begin{nr} \hfill 
$ \overline{R}(X,Y,Z,U)-\overline{R}(Z,U,X,Y)=<[\eta_X,\eta_Y]Z,U>-<[\eta_Z, \eta_U]X,Y>\hfill $
\end{nr}
whever $X,Y,Z,U$ are in $TM$. We end this section by recalling a technical consequence of proposition 2.1. 
\begin{lema}\cite{Nagy} Let $V,W$ belong to ${\cal{V}}$ and $X,Y$ be in $H$. We have : \\
(i) $\overline{R}(V,W)\eta=0$\\
(ii) $[\overline{R}(X,Y), \eta_V]=\eta_{\beta_V(X,Y)}$ where $\beta_V(X,Y)=\eta_{\eta_VY}X-\eta_{\eta_VX}Y$. 
\end{lema}
\section{Special ${\cal{AK}}_2$-manifolds}
This section will be devoted to devellop a number of preliminary results to be used in the proof of theorem 1.1. 
We begin by recalling the definition of 
special ${\cal{AK}}_2$-manifolds which is of esentially algebraic nature. 
\begin{defi}
Let $(M^{2n},g,J)$ be in the class ${\cal{AK}}_2$. It is said to be  special if and only if  $\eta_{{\cal{V}}} {\cal{V}}=H$ where 
$H$ is the K\"ahler nullity of $(g,J)$ and ${\cal{V}}$ its orthogonal complement in $TM$. 
\end{defi}
Algebraically speaking, the special condition ensures the vanishing of the torsion on ${\cal{V}}$, and also the symmetry of the restriction to ${\cal{V}}$ of the 
tensor $\eta$. We equally note that for a special ${\cal{AK}}_2$-manifold we also have : 
$$ (\nabla_{{\cal{V}}}J)H={\cal{V}}.$$
From a geometric viewpoint, the special condition tells us that the integral manifolds of the integrable distribution  ${\cal{V}}$ inherits from $(M,g,J)$ the structure of 
K\"ahler manifolds. \par
In the rest of this section we will work on a given ${\cal{AK}}_2$-manifold $(M^{2n},g,J)$. The notations of the previous definition 
are to used without further comment. We are going to investigate the action of the curvature tensor $\overline{R}$ on the decomposition $TM={\cal{V}} \oplus H$. 
Our starting point is the following intermediary result. 
\begin{lema}
Let $(M^{2m},g,J)$ be a special ${\cal{AK}}_2$ manifold. Then the following holds :
\begin{nr} \hfill 
$ 2 \overline{R}(V_3, V_4, V_2, \eta_{V_1}X)=\overline{R}(V_3,V_4, X,
\eta_{V_1}V_2)\hfill $
\end{nr}
for all $V_i, 1 \le i \le 4$ in ${\cal{V}}$ and $X$ in $H$.
\end{lema}
{\bf{Proof}} : \\
We will make use of the following formula from \cite{Nagy}, a consequence of the second Bianchi identity for 
the connection $\overline{\nabla}$: 
\begin{nr} \hfill 
$\overline{R}(\eta_{V_2}X,V_1,V_2,V_3,V_4)-
\overline{R}(\eta_{V_1}X, V_2, V_3, V_4)=-<[\eta_{V_3}, \eta_{V_4}]X,T_{V_1}V_2> \hfill $ 
\end{nr}
whenever $V_i, 1 \le i \le 4$ are in ${\cal{V}}$ and $X$ is in $H$. 
First, we note that under the special condition the right hand side of 
(3.2) vanishes so that we have $\overline{R}(\eta_{V_2}X,V_1,V_2,V_3,V_4)=
\overline{R}(\eta_{V_1}X, V_2, V_3, V_4)$. Now, using the symmetry property (2.2) we obtain : 
$$ \begin{array}{lr}
\overline{R}(V_3, V_4, \eta_{V_2}X, V_1)+<[\eta_{\eta_{V_2}X}, \eta_{V_1}]V_3,V_4>-
<[\eta_{V_3}, \eta_{V_4}]\eta_{V_2}X, V_1>= \\
\overline{R}(V_3, V_4, \eta_{V_1}X, V_2)+<[\eta_{\eta_{V_1}X}, \eta_{V_2}]V_3,V_4>-
<[\eta_{V_3}, \eta_{V_4}]\eta_{V_1}X, V_2>
\end{array} $$
Using the vanishing of the torsion $T$ on ${\cal{V}}$, a standard verification leads to 
$ <[\eta_{\eta_{V_2}X}, \eta_{V_1}]V_3,V_4>-
<[\eta_{V_3}, \eta_{V_4}]\eta_{V_2}X, V_1>=0$ hence 
\begin{nr} \hfill
$\overline{R}(V_3, V_4, \eta_{V_2}X, V_1)=\overline{R}(V_3, V_4, \eta_{V_1}X, V_2).  \hfill $
\end{nr}
But $(\overline{R}(V_3,V_4).\eta)(V_2,X)=0$ for $i=1,2$ (see lemma 2.1, (i)). Plugging this in the previous equation 
we obtain : 
$$ \begin{array}{lr}
<\eta_{\overline{R}(V_3,V_4)V_2}X,V_1>+<\eta_{V_2}\overline{R}(V_3,V_4)X, V_1>=
\overline{R}(V_3, V_4, \eta_{V_1}X, V_2).
\end{array}$$
Invoking proposition 2.2, (iii), one finds that the operator $\overline{R}(V_3,V_3)$ preserves ${\cal{V}}$ and $H$. Using again 
the vanishing of the torsion on ${\cal{V}}$ we have 
$$\begin{array}{lr}
<\eta_{\overline{R}(V_3,V_4)V_2}X,V_1>=-<X, \eta_{\overline{R}(V_3,V_4)V_2}V_1>=\\
-<X,  \eta_{V_1}\overline{R}(V_3,V_4)V_2>=\overline{R}(V_3,V_4, V_2, \eta_{V_1}X)
\end{array}$$
and the conclusion is now immediate
$\blacksquare$ \\ \par
The relation in lemma 3.1 shows that the restriction of $\overline{R}$ to ${\cal{V}}$ is completely determined 
by the mixed curvature terms of type $\overline{R}(V,W,X,Y)$ with $V,W$ in ${\cal{V}}$ and $X,Y$ in $H$. To 
investigate these terms we introduce now the 
{\bf{configuration}} tensor $A : H \times H \to {\cal{V}}$ by setting : 
$$ \overline{\nabla}_XY=\tilde{\nabla}_XY+A_XY$$ 
for all $X,Y$ in $H$. In a similar way, we define $B : H \times {\cal{V}} \to {\cal{V}}$ by 
$$ \overline{\nabla}_XV=\tilde{\nabla}_XV+B_XV.$$
Since $H$ is integrable $A$ is a symmetric tensor, that is $A_XY=A_YX$ for all $X,Y$ in $H$. It is immediate to establish that 
$ [A_X, J]=0$ for all $X$ in $H$. Now, the parallelism of $\eta$ in the direction of ${\cal{V}}$ together with the caracterization in proposition 2.1 translates into 
additional algebraic properties of the tensor $A$, as follows. 
\begin{lema}Let $X,Y$ be in $H$ and $V,W$ in ${\cal{V}}$ respectively. We have : \\
(i) $B_X(\eta_VY)=\eta_V(A_XY)$. \\
(ii) $A_X(\eta_VW)=\eta_V(B_XW).$
\end{lema}
{\bf{Proof}} : \\
Having in mind proposition 2.2, (iii) it suffices to project on $H$ and ${\cal{V}}$ respectively the 
identities $(\overline{\nabla}_X \eta)(V,Y)=0$ and $(\overline{\nabla}_X \eta)(V,W)=0$.
$\blacksquare$ \\ \par
The configuration tensor $A$ can be used to compute parts of the curvature tensor $\overline{R}$ in the following way.
\begin{lema}
Let $V,W$ and $X,Y$ be vector fields in ${\cal{V}}$ and $H$ respectively. We have : \\ 
(i) $$\overline{R}(V,X,W,Y)=<W, (\overline{\nabla}_VA)(X,Y)>-<B_XV, B_YW> $$ 
(ii) $$ <W, (\overline{\nabla}_VA)(X,Y)>=<V, (\overline{\nabla}_WA)(X,Y)>.$$
(iii) $$ \overline{R}(V,W,X,Y)=-<B_XV, B_YW>+<B_XW, B_YV>$$
for all $V,W$ in ${\cal{V}}$ and $X,Y$ in $H$.
\end{lema}
{\bf{Proof}} : \\
The proof of (i) will be omitted since a standard computation using only the proposition 2.2, (iii). Now, (ii) comes from (i) by means of the symmetry property 
(2.2). To prove (iii) one uses the first Bianchi identity for $\overline{R}$ when noticing that the latter do not contains derivatives of the torsion, in virtue of proposition 
2.1 and proposition 2.2, (ii). 
$\blacksquare $ \\ \par
We will start now to compute parts of the curvature tensor $\overline{R}$. To begin with, let us define the symmetric, $J$-invariant, partial 
Ricci tensors $r_1 : {\cal{V}} \to {\cal{V}}$ and 
$r_2 : H \to H$ by setting : 
$$\begin{array}{lr}
\sum \limits_{ v_k \in {\cal{V}}}^{} \overline{R}(v_k, Jv_k)V=r_1(JV) \\ 
\sum \limits_{v_k \in {\cal{V}}}^{} \overline{R}(v_k, Jv_k)X=r_2(JX)
\end{array} $$
Then we have : 
\begin{pro} The partial Ricci tensors $r_1, r_2$ can be computed by the following formulas : \\
(i) $<r_1V,W>=\sum \limits_{ e_i  \in H}^{} <B_{e_i}V, B_{e_i}W>$ \\
(ii) $<r_2X,Y>=-2\sum \limits_{v_k \in {\cal{V}}}^{} <B_{X}v_k, B_{Y}v_k>$ \\
where $V,W$ are in ${\cal{V}}$ and $X,Y$ belong to $H$ and $\{ v_k \}, \{ e_i\}$ are arbitrary orthonomal basis in ${\cal{V}}$ and 
$H$ respectively.
\end{pro}
{\bf{Proof}} : \\
(ii) follows by a simple computation involving lemma 3.3, (iii). Let us prove (i). Using (3.3), lemma 3.1 and (i) we obtain that 
$$ \begin{array}{lr}
\sum \limits_{ v_k \in {\cal{V}}}^{} \overline{R}(v_k, Jv_k,V,\eta_{W}X)=\sum \limits_{ v_k \in {\cal{V}}}^{} \overline{R}(v_k, Jv_k,W,\eta_{V}X)=\frac{1}{2}
\sum \limits_{v_k \in {\cal{V}}}^{}\overline{R}(v_k,Jv_k, X, \eta_VW) \\
=-\sum \limits_{ v_k \in {\cal{V}}}^{}<B_Xv_k,B_{\eta_{V}W}(Jv_k)>. 
\end{array} $$
Now, we have :  
$$ \begin{array}{lr}
\sum \limits_{ v_k \in {\cal{V}}}^{}<B_Xv_k,B_{\eta_{V}W}(Jv_k)>=
-\sum \limits_{ v_k \in {\cal{V}}}^{} \sum \limits_{e_i \in H}^{}<e_i, B_{\eta_{V}W}(Jv_k)><v_k, 
A_{e_i}X> \vspace{2mm}\\
=\sum \limits_{e_i \in H}<A_{e_i}\eta_{V}W, JA_{e_i}X>.
\end{array} $$
Or using in an apropriate way lemma 3.2 we have 
$$ \begin{array}{lr}
<A_{e_i}\eta_{V}W, JA_{e_i}X>=<\eta_{V}B_{e_i}W, 
A_{e_i}JX>=-<B_{e_i}W, B_{e_i}\eta_{V}(JX)>=\\
<B_{e_i}W, B_{e_i}J \eta_{V}X>
\end{array}$$ 
and the conclusion is now straightforward. $\blacksquare$ $\\$ \par
We are now able to give the main technical result of this section. 
\begin{pro}
Let $(M^{2n},g,J), n \ge 2$ be a special ${\cal{AK}}_2$-manifold. Then the following holds : \\
\begin{nr} \hfill 
$ \Delta^{{\cal{V}}} \vert A \vert^2=-5\vert r_1 \vert^2-2\vert \overline{\nabla}_{{\cal{V}}}A \vert^2. \hfill $
\end{nr}
Here, $\overline{\nabla}_{{\cal{V}}}$ denotes the restriction of $\overline{\nabla}$ to ${\cal{V}}$ and $\Delta^{{\cal{V}}}$ is the corresponding partial 
Laplacian, acting on functions. 
\end{pro}
{\bf{Proof}} : \\
From lemma 3.3, (ii)  we deduce that 
\begin{nr} \hfill 
$(\overline{\nabla}_{JV}A)(JX,Y)=(\overline{\nabla}_VA)(X,Y) \hfill $
\end{nr} 
for all $V$ in ${\cal{V}}$ and $X,Y$ in $H$ respectively. We consider now the partial Laplacian $D^{{\cal{V}}}$, acting on $A$ by : 
$$ (D^{{\cal{V}}}A)(X,Y)=-\sum \limits_{v_k \in {\cal{V}}}^{} (\overline{\nabla}^2_{v_k,v_k}A)(X,Y)$$
for all $X,Y$ belonging to $H$, where $\{ v_k \}$ is an arbitrary local orthonormal basis of ${\cal{V}}$. Derivating (3.5) it follows that 
$$(D^{{\cal{V}}}A)(X,Y)=\frac{1}{2}J \sum \limits_{v_k \in {\cal{V}}}^{}(\overline{R}(v_k,Jv_k).A)(X,Y). $$
Using proposition 3.1, we obtain further : 
$$ (D^{{\cal{V}}}A)(X,Y)=-\frac{1}{2}(r_1A_XY+2A_{r_2X}Y+2A_Xr_2Y).$$
Taking the scalar product with $A$ gives now $<D^{{\cal{V}}}A,A>=-\frac{1}{2}(\vert r_1 \vert^2+4\vert r_2 \vert^2)$, or further 
$<D^{{\cal{V}}}A,A>=-\frac{5}{2}\vert r_1 \vert^2$, after noticing that $\vert r_1 \vert^2=\vert r_2 \vert^2$. Now , the standard Weitzenb\"ock formula 
gives 
$$ \frac{1}{2}\Delta^{{\cal{V}}} \vert A \vert^2=<D^{{\cal{V}}}A,A>-\vert \overline{\nabla}_{{\cal{V}}}A \vert^2$$ 
and the claimed formula follows now easily.
$\blacksquare$
\begin{rema}
(i) Another way of proving formula (3.4) is the following. Define an almost complex structure $I$ on $M$ by setting 
$$ I=J \ \mbox{on} \ {\cal{V}} \ \mbox{and} \ I=-J \ \mbox{on} \ H.$$ 
Then it can be shown that $(M^{2n},g,I)$ is almost K\"ahler, and one can even show after some calculation that $(g,I)$ belongs to the class 
${\cal{AK}}_3$. Then the use of Sekigawa's formula gives exactly (3.4). Of course, one has to use (3.1) and furthermore compute all the remaining 
curvature terms. Since the calculations are of more length we prefered the direct approach.\\
(ii) If the manifold $(M,g,I)$ belongs to the class ${\cal{AK}}_2$ then the function $\vert \nabla I \vert^2$ is known to be constant and 
by the previous proposition we get that $(g,I)$ is in fact a K\"ahler structure. \\
\end{rema}

\section{Proof of theorem 1.1}
In this section we will give the proof of the theorem 1.1. This will be done by 
showing that it is always possible to restrict, at least locally and in dimension at least $6$, the study of special ${\cal{AK}}_2$-
manifolds, to the case when the function $\vert A \vert^2$ is constant. This, together with proposition 3.2 of the 
previous section, will enable us to prove theorem 1.1. \par
As in the section 3 let $(M^{2n},g,J)$ be a special ${\cal{AK}}_2$-manifold, with K\"ahler nullity $H$ and let ${\cal{V}}$ be the distribution 
orthogonal to $H$. We will first study the integral manifolds of $H$. For every $X$ in $H$ define a linear map : 
$$ \gamma_X : {\cal{V}} \to {\cal{V}} \ \mbox{by} \ \gamma_XV=\eta_VX .$$
The vanishing of the torsion on ${\cal{V}}$ implies that $\gamma_X$ is symmetric for all $X$ in $H$. \par
Let us denote by $\tilde{R} $ the curvature tensor of the connection $\tilde{\nabla}$, where we recall that $\tilde{\nabla}$ is the orthogonal projection of 
$\overline{\nabla}$ onto the decomposition $TM={\cal{V}} \oplus H$ (see section 3). 
The maps $\gamma_X$ are in relation with the curvature of $H$ (with respect to the connection $\tilde{\nabla}$), 
as the following lemma shows.
\begin{lema} Let $X,Y,Z$ be in $H$ and $V,W$ in ${\cal{V}}$. We have : 
$$ \tilde{R}(X,Y, \eta_VW, Z)=<[[\gamma_X, \gamma_Y],\gamma_Z]V,W>.$$
\end{lema}
{\bf{Proof}} : \\
Using the definition of $\tilde{\nabla}$ we obtain after a short computation that 
\begin{nr} \hfill 
$ \overline{R}(X,Y,Z^{\prime},Z)=\tilde{R}(X,Y,Z^{\prime}, Z)+<A_YZ^{\prime}, A_XZ>-
<A_XZ^{\prime}, A_YZ>\hfill $
\end{nr}
for all $X,Y,Z, Z^{\prime}$ in $H$. Now, by lemma 2.1, (ii) one obtains 
\begin{nr} \hfill 
$ \overline{R}(X,Y, \eta_VW, Z)+\overline{R}(X,Y,W, \eta_VZ)=<\eta_{\beta_V(X,Y)}W,Z>. \hfill $
\end{nr}
But the symmetry property (2.2) yields to 
$$ \overline{R}(X,Y, W, \eta_VZ)-\overline{R}(W,\eta_VZ, X,Y)=-<[\eta_W, \eta_{\eta_VZ}]X,Y>$$
since $H$ is the K\"ahler nullity of $(g,J)$. The use of lemma 3.3, (iii) gives then 
$$ \begin{array}{rr}
\overline{R}(W, \eta_VZ, X, Y)= -<B_XW, B_Y(\eta_VZ)>+<B_X(\eta_VZ), B_YW>\\
=-<B_XW, \eta_V(A_YZ)>+<\eta_V(A_XZ), B_YW>\\
= <A_X(\eta_VW), A_YZ>-<A_XZ, A_Y(\eta_VW)>.
\end{array} $$
where we used succesivelly lemma 3.2, (i) and (ii). It follows that 
$$\overline{R}(X,Y,W, \eta_VZ)=<A_X(\eta_VW), A_YZ>-<A_XZ, A_Y(\eta_VW)>-<[\eta_W, \eta_{\eta_VZ}]X,Y>
.$$ Using this in (4.2) and taking $Z^{\prime}=\eta_VW$ in (4.1) we get 
$$\tilde{R}(X,Y, \eta_VW, Z)=<\eta_{\beta_V(X,Y)}W,Z>+<[\eta_W, \eta_{\eta_VZ}]X,Y>.$$ It remains to 
take into account, in the last equation, the definition of the maps $\gamma_U, U$ in $H$ and our lemma follows 
routineously.
$\blacksquare$ \par
Our basic tool in the study of special ${\cal{AK}}_2$-manifolds will be the following intermediary result showing that 
a partial holonomy reduction of $H$ extends in a canonical way to a holonomy reduction over $TM$.
\begin{pro}
Let $(M,g,J)$ be a special ${\cal{AK}}_2$-manifold. Assume that we have an orthogonal, $J$-invariant decomposition 
$H_1 \oplus H_2$ which is also $\tilde{\nabla}$-parallel (inside $H$). Then : \\
(i) If we put ${\cal{V}}_i=\eta_{{\cal{V}}} H_i$ for $i=1,2$ then we have an orthogonal and $J$-invariant 
decomposition 
$${\cal{V}}_1 \oplus {\cal{V}}_2={\cal{V}}.$$
Moreover we have $\eta_{{\cal{V}}_1} {\cal{V}}_2=0$ and $\eta_{{\cal{V}}_i} {\cal{V}}_i=H_i, i=1,2$.\\
(ii) The decomposition 
$$ TM=({\cal{V}}_1 \oplus H_1) \oplus ({\cal{V}}_2 \oplus H_2)$$
defines a local splitting of $M$ into the Riemannian product of two special ${\cal{AK}}_2$-manifolds, with corresponding K\"ahler nullities 
$H_1$ and $H_2$.
\end{pro}
{\bf{Proof}} : \\
(i) Let $X_1$ and $X_2$ be in $H_1$ and $H_2$ respectively. Then the partial parallelism of $H_1$, together with the symmetry 
property of $\tilde{R}$ (a consequence of (4.1) and (2.2)) ensures that 
$\tilde{R}(X_1,X_2, \eta_VW, Z)=$ for all $V,W$ in ${\cal{V}}$ and 
$Z$ in $H$. Then, by the previous lemma we obtain 
$$ [[\gamma_{X_1}, \gamma_{X_2}],\gamma_Z]=0$$
for all $Z$ in $H$. Taking $Z=X_1$ we find that 
$$ \gamma_{X_1}^2 \gamma_{X_2}+\gamma_{X_2}\gamma_{X_1}^2=
2 \gamma_{X_1} \gamma_{X_2} \gamma_{X_1}.$$
We change now $X_2$ in $JX_2$ in the previous equation and take into account that 
$\gamma_{JX}=\gamma_XJ=-J\gamma_X$. It follows that 
$$ \gamma_{X_1}^2 \gamma_{X_2}+\gamma_{X_2}\gamma_{X_1}^2=
-2 \gamma_{X_1} \gamma_{X_2} \gamma_{X_1}$$
hence we must have 
$$\gamma_{X_1}^2 \gamma_{X_2}+\gamma_{X_2}\gamma_{X_1}^2=
 \gamma_{X_1} \gamma_{X_2} \gamma_{X_1}=0.$$
This yields to $\gamma_{X_1}^3 \gamma_{X_2}=0$ and since $\gamma_{X}$ is a symmetric operator 
for all $X$ in $H$ we get that $\gamma_{X_1} \gamma_{X_2}=0$. But this fact is easily seen to be equivalent to the orthogonality 
of the spaces ${\cal{V}}_1=\eta_{{\cal{V}}}H_1$ and 
${\cal{V}}_2=\eta_{{\cal{V}}}H_2$. The remaining afirmations of (i) are direct consequences of this fact.\\ 
(ii) We are going to prove first that the distribution ${\cal{V}}_1 \oplus H_1$ is $\overline{\nabla}$-parallel. Let 
$U$ be in $TM$ and $V,W$ in ${\cal{V}}_1$. As $(\overline{\nabla}_U \eta)(V,W)=(\overline{\nabla}_V\eta)(U,W)=0$ we 
get that 
$$ \overline{\nabla}_U(\eta_VW)=\eta_{\overline{\nabla}_UV}W+\eta_V(\overline{\nabla}_UW)$$ 
belongs to $\eta_{TM}{\cal{V}}_1+\eta_{{\cal{V}}_1}TM={\cal{V}}_1 \oplus H_1$. Since $H_1=\eta_{{\cal{V}}_1}{\cal{V}}_1$ 
we conclude that $\overline{\nabla}_UX$ belongs to ${\cal{V}}_1 \oplus H_1$ for all $X$ in $H_1$. 
Take now $V$ in ${\cal{V}}_1$ and $X$ in $H_1$. As before, we have : 
$$  \overline{\nabla}_U(\eta_VX)=\eta_{\overline{\nabla}_UV}X+\eta_V(\overline{\nabla}_UX)$$ 
belongs to $\eta_{TM}H_1+\eta_{{\cal{V}}_1}TM={\cal{V}}_1 \oplus H_1$ and using the fact that 
$\eta_{{\cal{V}}_1}H_1={\cal{V}}_1$ we conclude that ${\cal{V}}_1 \oplus H_1$ is $\overline{\nabla}$-parallel. In the same 
way it can be proven that ${\cal{V}}_2 \oplus H_2$ is $\overline{\nabla}$-parallel. Now, if $E_i={\cal{V}}_i \oplus H_i$ 
for $i=1,2$ it follows from (i) that $\eta_{E_i}E_i \subseteq E_i$ and $\eta_{E_i}E_j=0$ if $i \neq j$. This shows that 
$E_1$ and $E_2$ are in fact $\nabla$-parallel and the proof is finished.
$\blacksquare$ \\ \par
We will now study properties of the tensor $r_2 : H \to H$ defined in the previous section. 
\begin{lema} Let $(M^{2n},g,J)$ be a special ${\cal{AK}}_2$-manifold. Then : 
$$ (\tilde{\nabla}_Xr_2)Y=0$$
for all $X,Y$ in $H$.
\end{lema}
{\bf{Proof}} : \\
We will make use of the following formula which has been proven in \cite{Nagy} : 
\begin{nr} \hfill
$ (\overline{\nabla}_X \overline{R})(V_1,V_2,V_3,V_4)=0 \hfill $
\end{nr}
 for all $X$ in $H$ and $V_i $ in ${\cal{V}}, 1 \le i \le 4$. It follows easily that 
$(\tilde{\nabla}_X \overline{R})(V_1,V_2,V_3,V_4)=0$. Now, we recall that lemma 3.1 states that 
$$ 2 \overline{R}(V_3,V_4, V_2, \eta_{V_1}Y)=\overline{R}(V_3,V_4,Y, \eta_{V_1}V_2)$$ 
whenever $Y$ belongs to $H$ and $V_i$  in ${\cal{V}}, 1 \le i \le 4$. To conclude, it suffices to derive the last equation in the direction 
of $X$ in $H$, with 
respect to the connection $\tilde{\nabla}$, and next take (4.3) into account. 
$\blacksquare$ \\ \par
We will give now the proof of the theorem 1.1 stated in the introduction. \\
$\\$
{\bf{Proof of theorem 1.1}} \\
Let $U$ the open dense set of $M$ where we have an orthogonal splitting : 
$$ H=H_1 \oplus \ldots H_p $$
where $H_i, 1 \le i \le p$ are the eigenbundles of $r_2$ with corresponding eigenfunctions $\lambda_i, 1 \le i \le p$. Note 
that by lemma 4.2 we have in the standard way $X.\lambda_i=0$ for all $X$ in $H$ and $1 \le i \le p$ hence the distributions $H_i$ are 
$\tilde{\nabla}$-parallel, inside $H$. We set 
${\cal{V}}_i=\eta_{{\cal{V}}}H_i, 1 \le i \le p$ and use proposition 
4.1 to obtain an orthogonal, $J$-invariant and $\nabla$-parallel decomposition : 
\begin{nr} \hfill
$ TM=\bigoplus  \limits_{i=1}^{p} ({\cal{V}}_i \oplus H_i). \hfill $
\end{nr}
Of course, each factor corresponds to a special ${\cal{AK}}_2$-manifold such the corresponding tensor $r_2$ has the 
K\"ahler nullity as eigenspace. Therefore it remains us to consider this situation. \par
Suppose that $(M^{2n},g,J)$ is special ${\cal{AK}}_2$-manifold, with decomposition $TM={\cal{V}} \oplus H$ and such 
that $r_2=\lambda \cdot 1_H$. Then using (3.1) and the fact that $\eta_{{\cal{V}}}H={\cal{V}}$ we obtain that 
$$ \sum \limits_{v_k \in {\cal{V}}}^{} \overline{R}(v_k, Jv_k, V,W)=-\frac{\lambda}{2}<JV,W>$$ 
for all $V,W$ in ${\cal{V}}$. If $dim_{ \mathbb{R} }{\cal{V}}=2$ the special condition implies that $dim_{ \mathbb{R} }{H}=2$ and we know by the work in \cite{Apo1} that the torsion 
has to be parallel. If the dimension is greater, using the second Bianchi identity for $\overline{\nabla}$ on ${\cal{V}}$, exactly in 
the way one shows that a manifold of dimension greater than $3$, with Ricci tensor proportional to the metric 
tensor is Einstein, one obtains that $V.\lambda=0$ for all $V$ in ${\cal{V}}$. But using proposition 3.2 it follows 
immediately that the configuration tensor $A$ vanishes, ensuring the $\overline{\nabla}$-parallelism of the decomposition 
$TM={\cal{V}} \oplus H$. Moreover, by propositions 2.1 and 2.2, (ii) we obtain that the torsion is also $\overline{\nabla}$-parallel. \par
Now it clear that each factor of the decomposition (4.4) has parallel torsion and it follows that the torsion of 
$M$ is parallel over $U$ and by continuity over $M$, hence the first part of theorem 1.1 is proved. It remains us to show that $(M^{2n},g,J)$ is locally $3$-symmetric. The vanishing 
of $A$ implies by lemmas 3.1 and 3.3 that curvature terms of the form $\overline{R}(V_1, V_2, V_3, V_4), \overline{R}(V_1,V_2,X,Y)$ and $\overline{R}(V_1,X,V_2,Y)$ where 
$V_i, 1 \le i \le 4$ are in ${\cal{V}}$ and $X,Y$ in $H$ must be all equal to $0$. Furthermore, the restriction of $\overline{R}$ to $H$ is computed by lemma 4.1 and using the parallelism 
of the torsion it is an easy exercise to see that $\overline{\nabla} \overline{R}=0$, in other words $\overline{\nabla}$ is an Ambrose-Singer connection. We conclude now by 
\cite{Tri}. 
$\blacksquare$ \\ \par 
It is a good place now to introduce, in view of future use, the following definition.
\begin{defi}
A locally $3$-symmetric space of type $I$ is a special ${\cal{AK}}_2$ with parallel torsion.
\end{defi}
\begin{rema}
From the proof of theorem 1.1 we get the explicit dependence on the torsion of the curvature tensor $\overline{R}$ of a space of type I. This can be used to get 
algebraic caracterizations as a homogeoneous space of such a manifold. Since this discussion is beyond the scope of the present paper it will be omitted. 
\end{rema}
$\\$
{\bf{Proof of theorem 1.3}} : \\
It is easy to see that an almost K\"ahler structure satisfying the conditions in (i) or (ii) has to satisfy the second Gray condition and therefore must have parallel torsion. We conclude 
by recalling (cf. \cite{Nagy}) that any Einstein manifold supporting a compatible almost K\"ahler structure with parallel torsion is K\"ahler.

\begin{flushright}
Paul-Andi Nagy \\
Institut de Math\'ematiques \\
rue E. Argand 11, CH-2007, Neuch\^atel \\ 
email : Paul.Nagy@unine.ch
\end{flushright}
\end{document}